# THE TOPOGRAPHY OF MULTIVARIATE NORMAL MIXTURES[1]

By Surajit Ray and Bruce G. Lindsay

*University of North Carolina at Chapel Hill
and Pennsylvania State University*

Multivariate normal mixtures provide a flexible method of fitting high-dimensional data. It is shown that their topography, in the sense of their key features as a density, can be analyzed rigorously in lower dimensions by use of a ridgeline manifold that contains all critical points, as well as the ridges of the density. A plot of the elevations on the ridgeline shows the key features of the mixed density. In addition, by use of the ridgeline, we uncover a function that determines the number of modes of the mixed density when there are two components being mixed. A followup analysis then gives a curvature function that can be used to prove a set of modality theorems.

## 1. Introduction.

1.1. *The topography of a density.* Fitting a mixture model offers a primary data reduction through the number, location and shape of its components. However, in more complex settings we would like to know more about how the components interact to describe an overall pattern of density. What, for example, is the modal structure, or in a richer sense, the configuration of major features? The goal of this paper is to develop new insights into the topography of multivariate normal mixture densities, with the special aim of providing tools that are useful even in high-dimensional data.

Description of a multimodal density is challenging even in one dimension. For unimodal models, the density shape might be described through concepts like skewness and kurtosis. When the density is multimodal, the emphasis is usually shifted to the number and location of modes, since

---

Received April 2004; revised November 2004.

[1]Supported by NSF Grant DMS-01-04443.

---

the modes are the dominant feature, and are themselves potentially symptomatic of underlying population structures. A common approach to multimodal data structures is to use a mixture model because it provides a decomposition of the sampled population into a set of homogeneous components in a way that is consistent with the multimodal density configuration.

However, the set of modes is not in one-to-one correspondence with the distinct components. For example, in a univariate normal mixture model, two components can be similar enough that their mixture creates a single mode. If the population is well described by such a two component but-unimodal density, the analyst can construct two competing hypotheses. One is that the population has two homogeneous groups with normal shape, and the other is that there is just one group and it has a unimodal but nonnormal shape. Knowledge of these competing hypotheses could then lead to further scientific investigation.

Our purpose here is to describe the complex density shapes that can arise in a multivariate data set. To do so, we will appeal to the language and imagery of the earth's topography. Suppose we wish to describe the main features of a contour plot of a bivariate density $f(x, y)$. We equate this to the problem of describing the surface features of a land mass, where the elevation at a point $(x, y)$ is equated with the bivariate density $f(x, y)$. The local maxima of the density are the peaks, and their location, together with elevation, provides a first-order description of topography. But in a richer sense, mountains are usually aggregated into mountain ranges, in which the neighboring peaks are connected through ridges. The perceived separation of two neighboring peaks is then determined by the elevation at the lowest point on this ridge, the saddlepoint between them.

Here we will show how to create such a description for a mixture of high-dimensional normals or similar distributions, like multivariate-$t$. Our results are closely related to ideas in Morse theory, which is the mathematical study of the topology of surfaces based upon their differential structure [15, 16]. We will point out these relationships along the way.

1.2. *Relevant literature.* The literature on determination of the number of modes in normal mixture models has focused primarily on univariate mixtures. In fact, there is a simple description of modality when one is mixing two univariate components. de Helguero [5] determined necessary and sufficient conditions for bimodality in the mixture of two univariate normals with equal variances and mixing proportions. Later, conditions for bimodality in the mixture of univariate normal distributions with unequal variance and unequal mixing proportions were studied by Eisenberger [6], Behboodian [1] and Robertson and Fryer [20]. Kakiuchi [10] and Kemperman [11] addressed conditions for bimodality using nonnormal component densities.

Moreover, for the mixture of univariate normals, Robertson and Fryer [20] developed explicit formulae that one can use to determine if there are one or two modes, and further showed that there can be no more than two. This work, which can be shown to be a corollary of one of our results (Theorem 3), shows that modality is determined not only by the separation between the normal components, but also by the mixing proportion $\pi$, with values near 0 and 1 tending to suppress extra modes.

Some approaches to understand the modal behavior of multivariate normal mixtures can be found in recent machine learning literature (see [3]), but their results are limited to all components having the same covariance structure.

1.3. *Our results.* At this stage let us develop our notation. We will denote the dimension of the multivariate density by $D$ and the number of components of the mixture by $K$. All our topographical methods work for arbitrary $D$, but, for dimensionality reasons, our primary focus will be on cases where the number of components is $K = 2$ or 3. However, we also obtain a number of results for arbitrary $K$. As long as $K - 1$ is less than $D$, we can provide a dimension reduction to the problem.

In particular, if $K = 2$, a tremendous dimension reduction is possible, from $D$ (arbitrary) down to *one.* In Section 2 we show that the problem of finding the modes for $K$ components can be reduced to examining the density values on a $(K - 1)$-dimensional manifold of $\Re^D$, which we call the *ridgeline surface.* Surprisingly, the *ridgeline surface* does not depend on the mixing proportion $\pi$. We then design graphical and analytical tools to describe the density on this manifold. One of our key results is that this surface has a ridge-like property that enables one to determine global features from those that are local to the ridgeline surface.

In Section 3 we demonstrate plots of the density on the ridgeline, which we call *elevation plots.* One of the surprising results in our investigation is that, in a mixture of two multivariate normals, the statement "a mixture of two normals cannot have more than two modes" is false, unlike the univariate problem.

Next, in Section 4, we go into a deeper analysis of the modal structure. For $K = 2$ we have constructed a function $\Pi(\cdot)$ that does not depend on $\pi$, but whose plot can be used to determine the number and location of modes for each value of $\pi$. In Section 5 we develop a geometric analysis which shows that there exists a fundamental curvature function $\kappa(\cdot)$ whose *zeroes* determine the modality potential of a pair of component densities. For $D = 1$ this can be used to prove the Robertson and Fryer [20] results. For larger $D$ this result can be used to prove modality results for certain important special cases, such as equal or proportional covariance matrices.

Finally, Section 6 provides interpretation of our results, possible ways of generalizing these results and lists some unanswered questions on the topography of multivariate normal mixtures.

**2. The ridgeline manifold.** A $K$-component mixture of $D$-dimensional normals can be represented by the probability density function

$$g(\mathbf{x}) = \pi_1 \phi(\mathbf{x}; \boldsymbol{\mu}_1, \Sigma_1) + \pi_2 \phi(\mathbf{x}; \boldsymbol{\mu}_2, \Sigma_2) + \cdots + \pi_K \phi(\mathbf{x}; \boldsymbol{\mu}_K, \Sigma_K),$$

(1)

$$\mathbf{x} \in \Re^D,$$

where $\pi_j$ is the mixing proportion of component $j$, $\pi_j \in [0,1]$, $\sum_{j=1}^K \pi_j = 1$, and $\phi(\mathbf{x}; \boldsymbol{\mu}, \Sigma)$ is the density of a multivariate normal distribution with mean $\boldsymbol{\mu}$ and variance $\Sigma$. We will sometimes use $\phi_j(\mathbf{x})$ as shorthand notation for $\phi(\mathbf{x}; \boldsymbol{\mu}_j, \Sigma_j)$, and call $\phi_j$ the $j$th component density.

Our goal in this section is to show that, for any $D$-dimensional, $K$-component normal mixture, we can define a $(K-1)$-dimensional surface, which is guaranteed to include all the critical points (modes, antimodes and saddlepoints) of the $D$-dimensional mixture density. We will then show that plotting the density values in this $(K-1)$-dimensional space is completely informative about the main topographic features of the density in $D$-dimensional space.

2.1. *The $(K-1)$-dimensional ridgeline manifold.* To find this manifold, we examine the structure of the critical points of $g$ (i.e., all values of $\mathbf{x}$ where the first derivative of $g$ is equal to 0). First we need to introduce some terminology:

DEFINITION 1. The $(K-1)$-dimensional set of points

$$\mathcal{S}_K = \left\{ \boldsymbol{\alpha} \in \Re^K : \alpha_i \in [0,1], \sum_{i=1}^K \alpha_i = 1 \right\}$$

(2)

will be called the *unit simplex*. The function $\mathbf{x}^*(\boldsymbol{\alpha})$ from $\mathcal{S}_K$ into $\Re^D$ defined by

$$\mathbf{x}^*(\boldsymbol{\alpha}) = [\alpha_1 \Sigma_1^{-1} + \alpha_2 \Sigma_2^{-1} + \cdots + \alpha_K \Sigma_K^{-1}]^{-1}$$

(3)

$$\times [\alpha_1 \Sigma_1^{-1} \boldsymbol{\mu}_1 + \alpha_2 \Sigma_2^{-1} \boldsymbol{\mu}_2 + \cdots + \alpha_K \Sigma_K^{-1} \boldsymbol{\mu}_K]$$

will be called the *ridgeline function*. It will sometimes be written as $\mathbf{x}^*_{\boldsymbol{\alpha}}$. The image of this map will be denoted by $\mathcal{M}$ and called the *ridgeline surface* or *manifold*. If $K = 2$, it will be called the *ridgeline*, as it is a one-dimensional curve.

For $K = 2$ the ridgeline can be represented as

$$(4) \qquad \mathbf{x}^*(\alpha) = [\bar{\alpha}\Sigma_1^{-1} + \alpha\Sigma_2^{-1}]^{-1}[\bar{\alpha}\Sigma_1^{-1}\boldsymbol{\mu}_1 + \alpha\Sigma_2^{-1}\boldsymbol{\mu}_2],$$

where $\alpha \in [0, 1]$ and $\bar{\alpha} = 1 - \alpha$. As $\alpha$ varies from 0 to 1, the image of the function $\mathbf{x}^*(\alpha)$ defines a curve from $\boldsymbol{\mu}_1$ to $\boldsymbol{\mu}_2$.

REMARK 1.  Note that the ridgeline function $\mathbf{x}^*(\boldsymbol{\alpha})$ and, hence, the manifold $\mathcal{M}$, depend on the means and variances of the component densities, but not the mixing proportions $\pi_j$.

We next show that all the modes and saddlepoints of the full density $g(\mathbf{x})$ must occur in $\mathcal{M}$. We will later show how their exact location depends on the values of the $\pi_j$'s.

THEOREM 1.  *Let $g(\mathbf{x})$ be the density of a $K$-component multivariate normal density as given by* (1). *Then all of $g(\mathbf{x})$'s critical values and, hence, modes, antimodes and saddlepoints, are points in $\mathcal{M}$.*

PROOF.  See the Appendix.  □

The ridgeline surface has a simple structure if the component variances $\Sigma_i$ are equal. The following result can be found in [3].

COROLLARY 1.  *If $\Sigma_i = \Sigma$, then the convex hull of the means $\boldsymbol{\mu}_j$ of the component densities contains all critical points of the density $g$.*

PROOF.  The manifold $\mathcal{M}$ equals the convex hull.  □

To illustrate the ridgeline curve, we consider the following simple example:

EXAMPLE 1 (Two components, three modes, unequal variance).  The mixture density with $D = 2$ and $K = 2$, and the parameters

$$\boldsymbol{\mu}_1 = \begin{pmatrix} 0 \\ 0 \end{pmatrix}, \qquad \Sigma_1 = \begin{pmatrix} 1 & 0 \\ 0 & 0.05 \end{pmatrix},$$

$$\boldsymbol{\mu}_2 = \begin{pmatrix} 1 \\ 1 \end{pmatrix}, \qquad \Sigma_2 = \begin{pmatrix} 0.05 & 0 \\ 0 & 1 \end{pmatrix}, \qquad \pi_1 = \pi_2 = \tfrac{1}{2}.$$

Figure 1 shows the contours of the density given in Example 1. Overlaid on the contour plot is the ridgeline curve, showing how it passes through the three modes and saddlepoints of the density $g$.

Remark 2. Morse theory is used in the analysis of terrains [4] and watersheds [17], albeit with varying terminology. One can define a "critical net" or "watershed" as a map that describes the terrain through the location of the critical points in its elevation function, together with one-dimensional curves called "separatrices" that connect these critical points. A flow line represents the path taken by water under gravity, and a separatrix is defined to be a flow line that connects two critical points. The flow line from a local maximum to a saddlepoint therefore creates a division of the terrain into two "catchment basins," as flow lines do not cross. Mathematically, a flow line is a line of steepest descent (or ascent when reversed), so that the flow line is always moving in the direction of the gradient. It is therefore always moving orthogonally to the elevation contours. For $K = 2$ our ridgeline curve has properties similar to the separatrices, as we shall show, but is not one itself. In particular, the separatrices for a mixture of two normals depend on the value of $\pi$, but the ridgeline curve does not.

2.2. *The ridgeline elevation plot.* The next step in our analysis is to consider the diagnostic properties of the *elevation plot*, which is a plot of the *ridgeline elevation function* defined by

$$h(\boldsymbol{\alpha}) = g(\mathbf{x}^*(\boldsymbol{\alpha})).$$

We start by considering the case where $K = 2$, so $\alpha$ is one-dimensional. Figure 2(a) shows the elevation plot of the distribution of Example 1. One might hope that the number and location of the modes and antimodes of the elevation plot might indicate to us the number and location of the modes and saddlepoints of the original density. We will make this rigorous through further analysis in the next subsection.

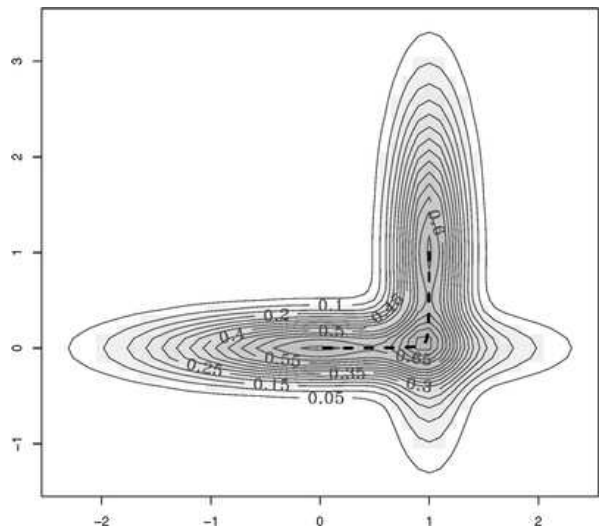

Fig. 1. *Contour plot and ridgeline curve (- - - - -) for the mixture density given in Example 1.*

REMARK 3. The elevation plot carries an inherent visual distortion relative to the original density plot. This distortion arises because the distance between two $\alpha$ values, say $\alpha_1$ and $\alpha_2$, may not accurately reflect the distances between $\mathbf{x}^*(\alpha_1)$ and $\mathbf{x}^*(\alpha_2)$ as measured along the ridgeline. For example, the saddlepoints in the contour plot of Figure 1 have a distance from the endpoints that is not well represented on the $\alpha$-scale in Figure 2(a). For $K = 2$ this distortion can be corrected by replacing the $\alpha$ scale with $L(\alpha)$, the arclength of the ridgeline path $\mathbf{x}^*(\alpha)$ from 0 to $\alpha$ [see Figure 2(b)].

REMARK 4. Theorem 1 can be easily generalized to other families of multivariate distributions. If a density characterized by parameters $(\boldsymbol{\mu}, \Sigma)$ depends on $\mathbf{x}$ only as a decreasing function of the Mahalonobis distance $(\mathbf{x} - \boldsymbol{\mu})'\Sigma^{-1}(\mathbf{x} - \boldsymbol{\mu})$, then the ridgeline manifold is exactly the same as given here. An example of this type is the multivariate $t$ distribution (see [13, 18]).

2.2.1. *Ridgeline-like properties.* We now verify that the modes and antimodes in the ridgeline elevation plots, such as in Figure 2, necessarily correspond to modes and saddlepoints in the original density. Recall that in terrain analysis [4] a separatrix is the line connecting the highest points along a ridge and separating drainage basins from one another. We have already indicated that our ridgelines are not true separatrices. However, confirmation that the ridgeline elevation plot is fully informative about the main features of $g$ arrives by showing that the ridgeline surface has local properties similar to a separatrix. This section is devoted to making this idea mathematically precise.

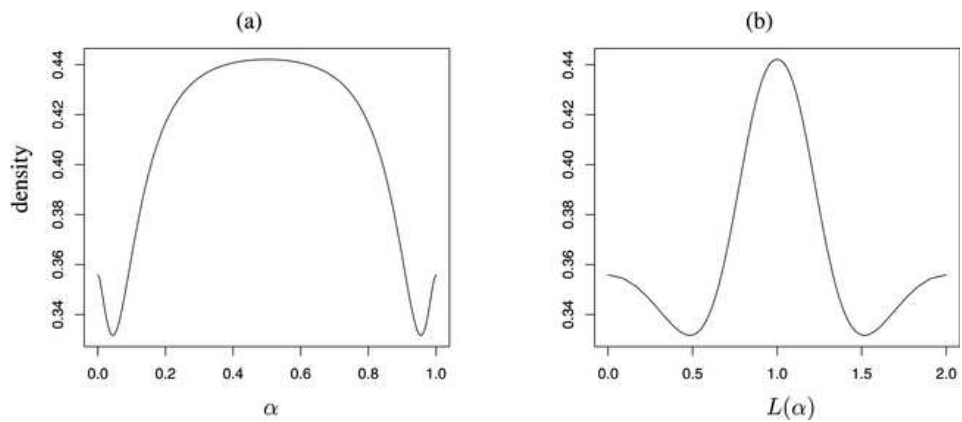

FIG. 2. *Ridgeline elevation for the bivariate normal mixture of Example 1 along the ridgeline path* $\mathbf{x}^*(\alpha)$, *expressed as a function of* (a) *parameter* $\alpha$ *and* (b) *the arclength* $L(\alpha)$. *Three local maxima representing the three modes of the density are visible near* $\alpha = 0$, 0.5 *and* 1 *in* (a), *which corresponds to* $L(\alpha) = 0$, 1 *and* 2 *in* (b).

If a person walks from peak A to peak B along a separatrix path, then this passage would be characterized by the fact that to the left and right, perpendicularly, the ridge falls away into two drainage basins. In other words, one stands at a local maximum in elevation of the mountain cross-section perpendicular to the path. Because of this, any time one reaches a point along the path that is locally maximal in elevation relative to the other points on the path, it must also be a local maximum to the entire nearby surface. By the same logic, when one reaches a local minimum along the path, it must be a saddlepoint to the surface. In this way one can infer the nature of the critical points on the whole surface from just the elevations along the path.

The motivation for our analysis arises from the following geometric property of the ridgeline curve $\mathbf{x}^*(\alpha)$ when $K = 2$. Consider any contour $\{\mathbf{x} : \phi_1(\mathbf{x}) = c\}$ of the first component density. This forms an ellipse. Provided that $\boldsymbol{\mu}_2$ is not inside this ellipse, one can create a contour of the second component, say $\{\mathbf{x} : \phi_2(\mathbf{x}) = d\}$, such that the two ellipses intersect at a single point $\mathbf{x}_0$. One can show that $\mathbf{x}_0$ is necessarily a point on the ridgeline curve, and, in fact, all points on the ridgeline curve are the "kissing points" of two such ellipses. Now, consider any point $\mathbf{x}$ that is on the hyperplane separating the two ellipses, but is not $\mathbf{x}_0$. Since it is in neither ellipse, it must have a smaller density value under $\phi_1$ and $\phi_2$ than $\mathbf{x}_0$ does. But this means that $g(\mathbf{x}) < g(\mathbf{x}_0)$, regardless of the value of $\pi$. That is, $\mathbf{x}_0$ is a local maximum relative to all points in the hyperplane. In fact, it must be a maximum relative to all points not in the two ellipses, independently of the value of $\pi$.

It is also clear (from its above geometric description) that the local direction of the ridgeline path never lies in the hyperplane. Hence, if derivative of the elevation along this path becomes zero at a point, there is a full rank set of directions in which the directional derivatives are zero and, hence, the gradient is zero at this point, making it a critical point.

To make this argument precise and to extend it to the multivariate case, we need to develop some notation. If $K = 2$, so that $\alpha$ is a scalar, the derivative vector $\mathbf{d}(\alpha) = d\mathbf{x}(\alpha)/d\alpha$ of the ridgeline points in the direction of the curve's travel. If $K > 2$, then $\boldsymbol{\alpha}$ is a vector, and the $K - 1$ derivative vectors of $\mathbf{x}(\boldsymbol{\alpha})$ with respect to $\alpha_1, \ldots, \alpha_{K-1}$, written $\mathbf{d}_1, \ldots, \mathbf{d}_{K-1}$, represent a basis for the space of possible directions of travel within the ridgeline surface. Using (3), these can be found as

$$\mathbf{d}_j = S_{\boldsymbol{\alpha}}^{-1}(\mathbf{v}_j - \mathbf{v}_K),$$

where for fixed $\alpha$ we have defined the matrix $S_{\boldsymbol{\alpha}} = \sum \alpha_j \Sigma_j^{-1}$ and vectors $\mathbf{v}_j = \Sigma_j^{-1}(\mathbf{x}^*(\boldsymbol{\alpha}) - \boldsymbol{\mu}_j)$. We next define a linear subspace of vectors that are orthogonal to the surface's direction vectors $\mathbf{d}_j$ in an appropriate sense,

$$\mathbf{W} = \{\mathbf{w} : \mathbf{w}' S_{\boldsymbol{\alpha}} \mathbf{d}_j = 0 \ \forall j = 1, \ldots, K-1\}.$$

THEOREM 2. *If $\mathbf{w} \in \mathbf{W}$, then along the path $\{\mathbf{x}(\boldsymbol{\alpha}) + \delta \mathbf{w} : \delta \in \Re\}$ the function $g(\mathbf{x})$ takes its maximum value at $\delta = 0$.*

PROOF. See the Appendix. □

COROLLARY 2. *Every critical point $\boldsymbol{\alpha}$ of $h(\boldsymbol{\alpha})$ corresponds to a critical point $\mathbf{x}(\boldsymbol{\alpha})$ of $g(\mathbf{x})$. A critical point $\boldsymbol{\alpha}$ of $h(\boldsymbol{\alpha})$ gives a local maximum of $g(\mathbf{x})$ if and only if it is a local maximum of $h(\boldsymbol{\alpha})$. If $D > K - 1$, so that the $h(\boldsymbol{\alpha})$ plot is a true dimension reduction, then $g(\mathbf{x})$ has no local minima, only saddlepoints and local maxima. In general, for $D \geq K - 1$, at a critical point of $h(\boldsymbol{\alpha})$ whose second derivative matrix has $m$ negative eigenvalues the function $g(\mathbf{x})$ will have a critical point whose second derivative matrix has an additional $D - K + 1$ negative eigenvalues corresponding to the dimension of the orthogonal directions $\mathbf{w}$.*

PROOF. The directional vectors, together with their orthogonal complement $\mathbf{W}$, span the space, and we know that the $\mathbf{W}$ vectors are all directions of local maximization. □

**3. Some illustrative examples.** Before we proceed further with the theory, we present some examples to illustrate the theory and methods developed to this point. We also use them to motivate the next set of theoretical developments.

3.1. *Elevation plots for $K = 2$.* Let us return to Figure 1. It appears from the contour shapes that there are three modes, all lying on the ridgeline. This shows that the multivariate normal case has a very different, and more complex, modal structure than the univariate. (Carreira-Perpiñán and Williams [3] also give an example in which there are three modes.) From Figure 2 we can see clearly that there are indeed three modes. (Note: the modes that appear to be at the endpoints actually occur a slight distance from the ends, although it is not visible in this resolution.) We also can see that the central mode is dominant, and the minima (corresponding to saddlepoints in Figure 1) are relatively shallow. Note that the contour plot in Figure 1 is not available unless $D = 2$, but the elevation plot works for any $D$.

This example above opens up a natural question regarding the existence of an upper bound to the number of modes in a two component mixture of multivariate normals. We have constructed another example where we can find four modes in three dimensions.

EXAMPLE 2 (Two components, four modes, unequal variance).   For $K = 2$ and $D = 3$, let the parameters be

$$\boldsymbol{\mu}_1 = \begin{pmatrix} 0 \\ 0 \\ 0 \end{pmatrix}, \qquad \Sigma_1 = \begin{pmatrix} 1 & 0 & 0 \\ 0 & 1 & 0 \\ 0 & 0 & 0.05 \end{pmatrix},$$

$$\boldsymbol{\mu}_2 = \begin{pmatrix} 1/\sqrt{2} \\ 2 \\ 1/\sqrt{2} \end{pmatrix}, \qquad \Sigma_2 = \begin{pmatrix} 0.05 & 0 & 0 \\ 0 & 1 & 0 \\ 0 & 0 & 1 \end{pmatrix}, \qquad \pi_1 = \pi_2 = \tfrac{1}{2}.$$

Figure 3 is the ridgeline elevation plot of the mixture in Example 2. It has four modes, although they are not clearly visible from the plot in Figure 3(a). We can more clearly see them in Figures 3(b) and 3(c), which show the ridgeline elevation plot for narrower ranges of values. Determination of the exact location in $\alpha$ of the modes can be done numerically using any nonlinear optimization software (we used the `nlm` function in R). In our case we found that the four modes are located at $\alpha = 0.00084, 0.137, 0.863$ and $0.99916$.

These examples raise the question: can we create further mathematical tools that will guide us in the determination of the number and location of modes? The answer is yes, but before tackling this we consider elevation plots in a higher dimension.

3.2. *Elevation plots for $K = 3$.*  For $K = 3$ the dimension of $\boldsymbol{\alpha}$ and, hence, the ridgeline surface, is two, so we suggest using a contour plot of $h(\boldsymbol{\alpha})$ in the $\boldsymbol{\alpha}$ coordinate system to look for critical points of the density.

To distinguish between the contour plot of an original density (possible only when $D \leq 2$) and its contour described along the ridgeline surface (available for arbitrary $D$), we will denote the former as the *density contour* and the latter as the *ridgeline contour*.

A distance-preserving way to represent the simplex $\mathcal{S}_3$ in $\Re^2$ is as an equilateral triangle. The vertices of the triangle correspond to the three equidistant points $\mathbf{e}_1 = (1, 0, 0)$, $\mathbf{e}_2 = (0, 1, 0)$ and $\mathbf{e}_3 = (0, 0, 1)$ of the simplex. Each point in the triangle corresponds to a point $(\alpha_1, \alpha_2, \alpha_3)$, where $\alpha_j$ equals the length of the perpendicular dropped to the side opposite to the vertex $\mathbf{e}_j$. We will use the symbol $\alpha_j$ at the corner $\mathbf{e}_j$ because the distance from the opposing baseline gives the $\alpha_j$ value. At the corner itself, $\alpha_j$ equals one.

We plot the ridgeline contours on this triangle. We have shown that the number of peaks on the contour plot is exactly the number of modes of the 3-component normal. We can also find exact positions of the modes by determining the simplicial coordinates $\boldsymbol{\alpha}$, then computing $\mathbf{x}^*(\boldsymbol{\alpha})$; these locations will depend on the values of $\pi$ that are used in $g$.

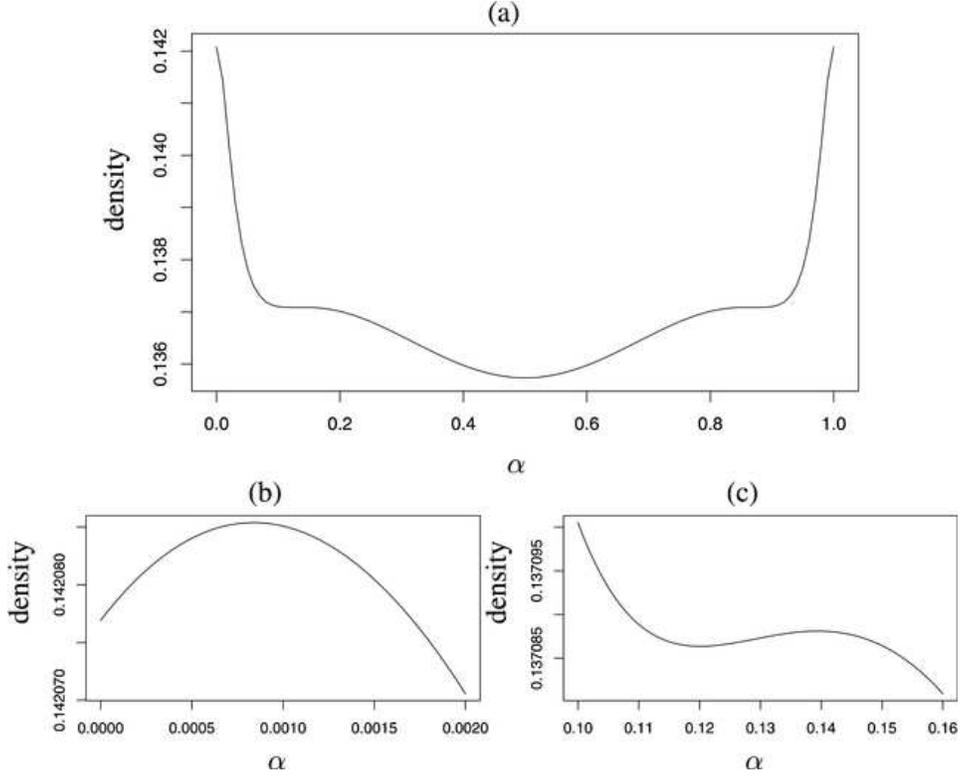

Fig. 3. *Ridgeline elevation plot of density of Example* 2. (a) *For the whole range of* $\alpha$, (b) *magnified for* $\alpha \in (0, 0.0002)$ *and* (c) *magnified for* $\alpha \in (0.1, 0.16)$.

EXAMPLE 3 (Three components, three modes, equal variance). For $K = 3$, $D = 2$, let the covariance matrix be common and the parameters be

$$\boldsymbol{\mu}_1 = \begin{pmatrix} 0 \\ 0 \end{pmatrix}, \qquad \boldsymbol{\mu}_2 = \begin{pmatrix} 0 \\ 3 \end{pmatrix}, \qquad \boldsymbol{\mu}_3 = \begin{pmatrix} 3 \\ 0 \end{pmatrix},$$

$$\Sigma = \begin{pmatrix} 1 & 0 \\ 0 & 1 \end{pmatrix}, \qquad \pi_1 = \pi_2 = \pi_3 = \tfrac{1}{3}.$$

Figure 4 shows the modality surface contour plot of Example 3. The three peaks along the three corners of the triangle are easily visible. In fact, these three modes lie very close to the corners, which implies that the three modes of the density $g$ are close to the three means of the component densities. Of course, for this example one could have done a contour plot of $g$ itself because $D = 2$.

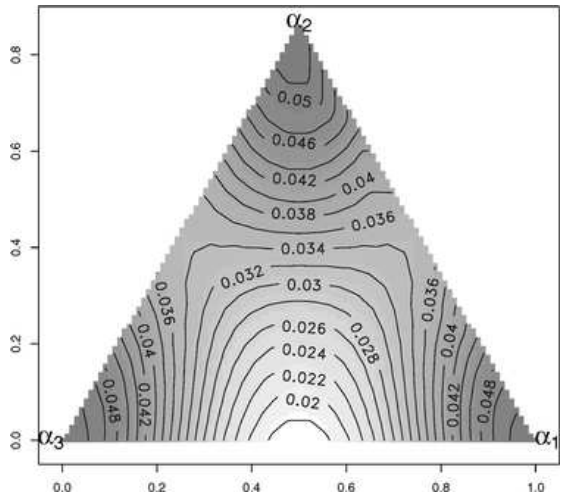

Fig. 4.  *Ridgeline contour plot of the density of Example* 3.

Now we move on to an unequal variance example.

EXAMPLE 4 (Three components, five modes, unequal variance).   For $K = 3$, $D = 2$, let the parameters be

$$\boldsymbol{\mu}_1 = \begin{pmatrix} 0 \\ 0 \end{pmatrix}, \qquad \boldsymbol{\mu}_2 = \begin{pmatrix} 1 \\ 1 \end{pmatrix}, \qquad \boldsymbol{\mu}_3 = \begin{pmatrix} 2 \\ 2 \end{pmatrix},$$

$$\Sigma_1 = \Sigma_3 = \begin{pmatrix} 1 & 0 \\ 0 & 0.05 \end{pmatrix}, \qquad \Sigma_2 = \begin{pmatrix} 0.05 & 0 \\ 0 & 1 \end{pmatrix}, \qquad \pi_1 = \pi_2 = \pi_3 = \tfrac{1}{3}.$$

If we carefully look at the ridgeline contour plot of Figure 5(a), we can find five modes, which can be verified from its density contour plot in Figure 5(b). The five modes are (i) near $\alpha_3 = 1$, (ii) near $\alpha_1 = 1$, (iii) near the centroid (one with height of the contour being 0.23), (iv) near $\alpha_1 = \alpha_2 = 0.5$ and (v) near $\alpha_2 = \alpha_3 = 0.5$. In this example, although Figures 5(a) and 5(b) carry the same modal information, we can see rather dramatically that the ridgeline contour plot distorts the distances and angles between the modes relative to the density plot.

EXAMPLE 5 (Iris data).   Next, we consider a three component mixture model fit to Fisher's iris data [7]. This is the dataset made famous by Fisher, who used it to illustrate principles of discriminant analysis. Data on four variables, namely, *Petal width*, *Petal length*, *Sepal width and Sepal length*, were collected on flowers of three iris species: Setosa, Verginica and Versicolor. Each species had 50 observations.

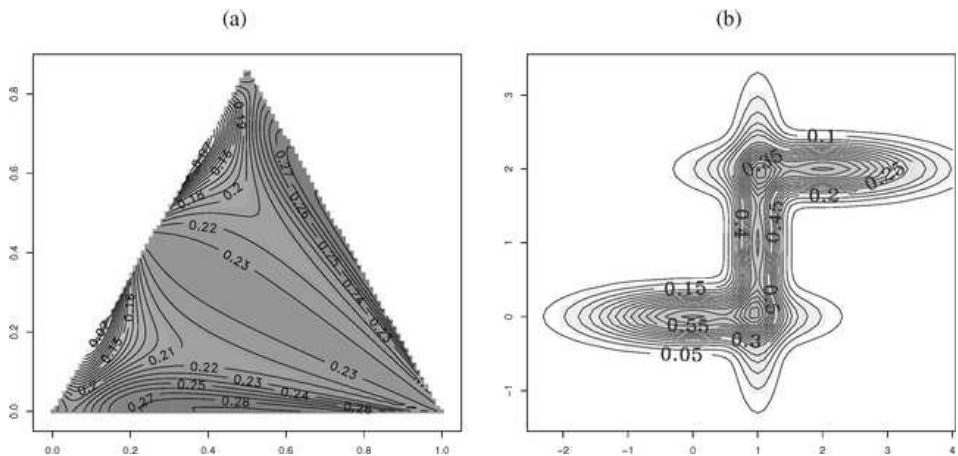

Fig. 5.  (a) *Ridgeline contour plot and* (b) *density contour plot of the three component mixture density of Example 4 with five modes.*

Since $D = 4$, direct contour plotting of $g$ is not available. The mixture models we used are the maximum likelihood fits to the data set assuming unequal variance. Examining the ridgeline contour plot in Figure 6(a), we conclude that the three component fit actually corresponds to three different modes, the modes being near the mean for each component, which corresponds to the three vertices in the simplex of Figure 6(a).

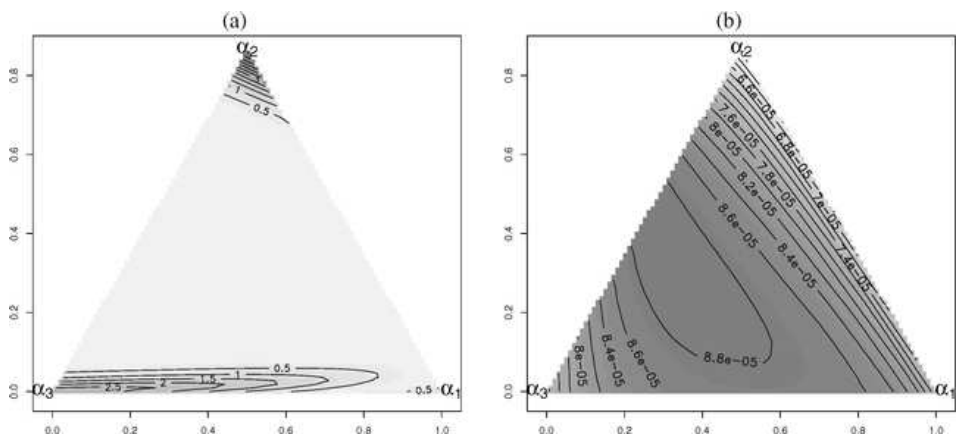

Fig. 6.  *Ridgeline contour plot of* (a) *Example 5 and* (b) *Example 6.*

EXAMPLE 6 (Egyptian skull data). This data consists of four measurements: *Maximal breadth*, *Basibregmatic height*, *Basialveolar length and Nasal height* of male Egyptian skulls from five different time periods (4000 BC, 3300 BC, 1850 BC, 200 BC, 150 AD). Thirty skulls were measured from each time period [21].

Here we analyze the three earliest time periods and fit a three component normal mixture with unequal variances. Examining the ridgeline contour plot [Figure 6(b)], we observe that the three components, pertaining to the three time periods, produce a single mode.

REMARK 5. More expressive detail for the two-dimensional plots of this section could have been obtained by displaying the critical net of the density using, for example, the approximation techniques of Danovaro et al. [4]. This would show the maxima, saddlepoints and separatrices over the manifold region based on evaluating the elevation at a finite network of points.

**4. The Π-plot.** Until this point we have focused on graphical techniques that are based on displaying the elevation of the density on the ridgeline. These techniques are quite elementary, and carry full information about the location and relative heights of the modes and saddlepoints. We now turn to a technique that focuses instead on the location of the modes and saddlepoints and not their elevations. By doing so, we can gain important insights into how the number and location of the modes depend on the values of $\pi$ for a given fixed set of component densities. For this section we will consider only the case $K = 2$. We will treat the component parameters as fixed throughout the analysis.

Recall the ridgeline curve $\mathbf{x}^*(\alpha)$ for $K = 2$, defined in (3). If $\mathbf{x}^*(\alpha)$ is a critical value of $h(\alpha)$, then it satisfies

$$h'(\alpha) = \pi \phi_1(\mathbf{x}^*(\alpha))' + \bar{\pi} \phi_2(\mathbf{x}^*(\alpha))' = 0,$$

where prime "′" denotes differentiation with respect to $\alpha$. Also, recall that, for fixed component parameters, the value of $\alpha$ completely specifies the vector $\mathbf{x}^*(\alpha)$ independently of $\pi$. For notational ease, any function of the $D$-dimensional vector $\mathbf{x}(\alpha)$ will be written as the same function of the scalar $\alpha$, for example, $\phi_1(\alpha)$ and $\phi_1(\mathbf{x}(\alpha))$ are the same. Solving the last displayed equation for $\pi$, and turning it into a function of $\alpha$, we get

$$\Pi(\alpha) = \frac{\phi_2'(\alpha)}{\phi_2'(\alpha) - \phi_1'(\alpha)}.$$

It follows that if $\alpha$ corresponds to a critical value for the density $g$, then it solves the "pi-equation"

$$(5) \qquad\qquad \Pi(\alpha) = \pi.$$

This gives us a recipe for finding the critical points of the density $g$, given a particular choice of $\pi$: first, we solve for the set of $\alpha$ that satisfy the pi-equation (5), and then, from these $\alpha$-solutions, calculate the critical points $\mathbf{x}^*(\alpha)$. If we have a unimodal density for a given $\pi$, there must be one and only one solution to the pi-equation, which corresponds to the mode of the density. A bimodal density will have three $\pi$-critical points, the first and the third (in the order of magnitude) corresponding to the two modes.

To aid in finding solutions to the pi-equation, we describe $\Pi(\alpha)$ on $\alpha \in [0, 1]$. We first claim that

$$\Pi(0) = 0, \qquad \Pi(1) = 1, \qquad \Pi(\alpha) \in [0, 1].$$

To show this, first notice that $\phi_1(\alpha)$ increases on the range $\alpha \in [0, 1)$ because $\phi_1'(\alpha) > 0$; also $\phi_1'(1) = 0$. On the other hand, $\phi_2(\alpha)$ is decreasing on range $\alpha \in (0, 1]$ as $\phi_2'(\alpha) < 0$ and $\phi_2'(0) = 0$. This establishes the result.

We can also derive the following simple calculation formula for $\Pi(\alpha)$:

$$(6) \qquad \frac{1}{\Pi(\alpha)} = \frac{\phi_2'(\alpha) - \phi_1'(\alpha)}{\phi_2'(\alpha)} = 1 + \frac{\alpha}{\bar{\alpha}} \frac{\phi_1(\alpha)}{\phi_2(\alpha)},$$

which can be verified by routine calculus.

As an example, let us examine the $\Pi$-plot of the two component bivariate normal mixture with three modes given in Example 1. As the mixing proportion in Example 1 is $\pi = 0.5$, we would draw a horizontal line across the $(\alpha, \Pi(\alpha))$ plot (Figure 7) at height $\pi = 0.5$. This line crosses the curve five times at ($\alpha = 0$(approx.), 0.004, 0.5, 0.996, 1(approx.)). Among these, $\alpha = 0$(approx.), 0.5, 1(approx.) correspond to the three modes, as was verified by the ridgeline elevation plot (see Figure 2).

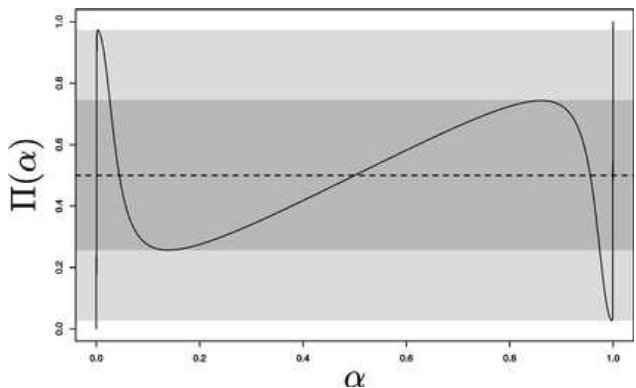

Fig. 7. $\Pi$-plot of Example 1. $- - - - -$ denotes mixing proportion $\pi = 0.5$. The ranges of $\pi$ for which the distribution has two and three modes are given by the light band and dark band, respectively.

The same Π-plot can also be used to determine the number of modes of $g$ as $\pi$ varies but the component densities remain unchanged. For example, if the means and variances are unchanged in Example 1, there will be three modes if and only if $0.256 < \pi < 0.744$ (dark band in Figure 7). Moreover, if $\pi < 0.0225$ or $\pi > 0.975$ (the unshaded region), then we will have a unimodal density; for all other $\pi$'s, we will have a bimodal density (light band).

Similarly, for Example 2 we can find the range of mixing proportions for which the distribution will have one, two, three or four modes by examining the plots in Figures 8(a) and 8(b). There is a narrow range, $\pi \in (0.49974, 0.50026)$, for which the parameters generate a distribution with four modes. The number of modes is at least two if $\pi \in (0.25, 0.75)$ and at least three if $\pi \in (0.489, 0.511)$.

When $K > 2$, Π-plots can be used to examine the modal structure of each pair of component densities, where the mixing proportions in the paired mixture is determined by the relative weight in the whole mixture. To illustrate, we re-examine (i) Example 4, (ii) Example 5 and (iii) Example 6. The ridgeline contour plots of the above examples were already presented in Section 3. Now we present the pairwise Π-plots of the three components fit for each dataset (Figure 9).

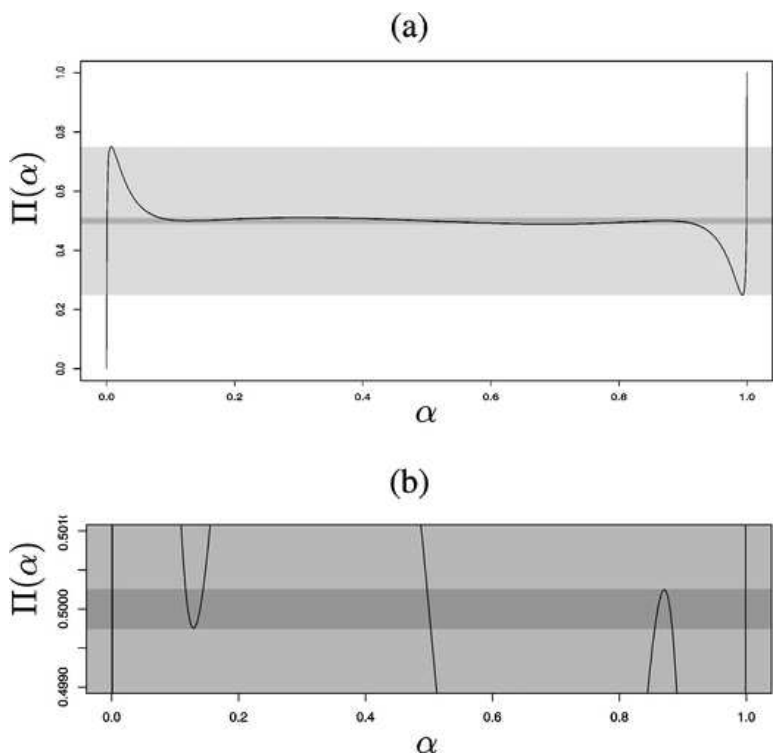

FIG. 8.  Π*-plot of Example* 2. *The range of* $\pi$ *for which the mixing distribution has two modes is given by the light band; with three modes it is given by the medium band; with four modes it is given by the dark band. The first two bands are visible from sub-plot* (a), *whereas the two darker bands are more clearly visible in sub-plot* (b), *which is magnified for* $\pi \in (0.499, 0.501)$.

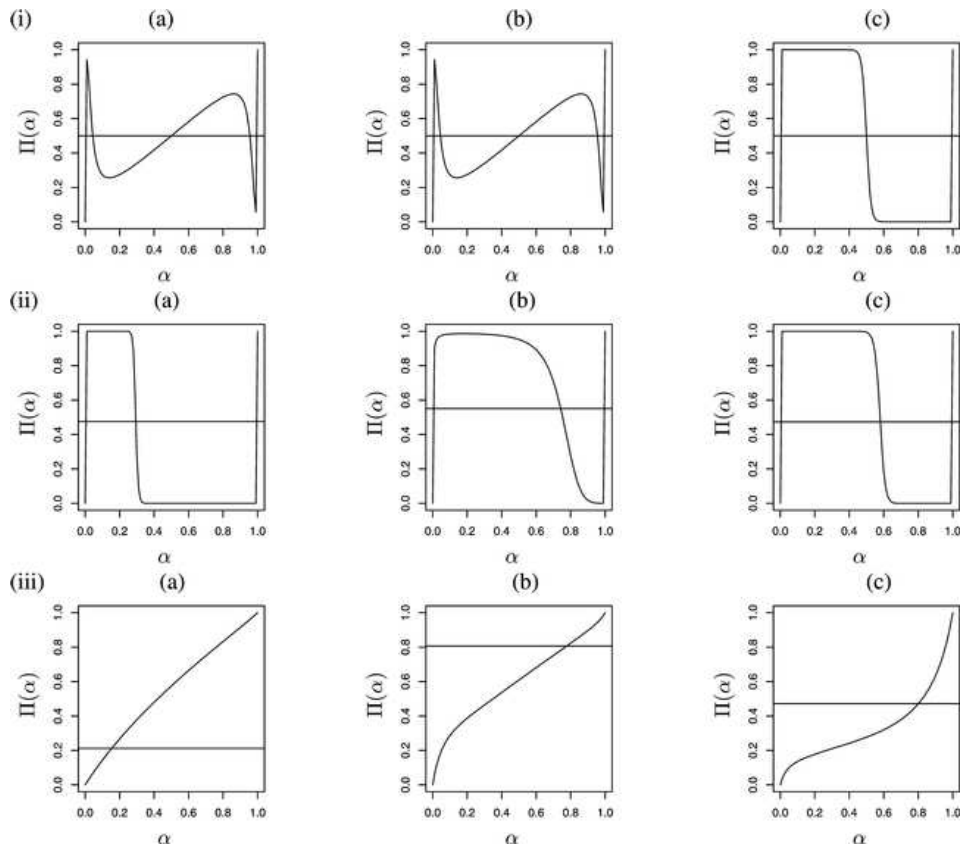

Fig. 9. Π-plots of pairs of component densities in the three-component mixture densities of (i) *Example 4*, (ii) *Example 5 and* (iii) *Example 6, for component pairs* (a) {1, 2}, (b) {2, 3}, (c) {3, 1}, *respectively, along with the appropriate mixing proportion represented by the horizontal line in each plot.*

For the pairwise Π-plots of Example 4 [Figure 9(i)], we observe that, while comparing pairs {1, 2} and {2, 3}, the horizontal line crosses the Π-curve five times, implying the presence of three modes. Similarly, examining the pairwise plots for the iris data [Figure 9(ii)], we see that each of the pairwise components is bimodal. Moreover, the plot tells us that the three components are well separated from each other for almost every set of mixing proportions. But, looking at the plots for the Egyptian skull data [Figure 9(iii)], we can easily note that all the plots imply that the pairwise components exhibit a single mode.

## 5. Analytic tools for detecting modality.

5.1. *The curvature function.* We have seen that the $\Pi$ function determines the modality structure of the mixture model as it depends on $\pi$. Our next step is to look more deeply into the properties of this function. Now if $\Pi$ is strictly increasing in $\alpha$, there will be exactly one solution to the pi-equation, and so the mixture densities are unimodal across all values of $\pi$. Multiple solutions in $\alpha$ of the pi-equation can only arise if $\Pi(\alpha)$ oscillates up and down as $\alpha$ increases. In particular, a careful analysis shows that $\Pi$ is increasing at a mode of $h(\alpha)$, and decreasing where $h(\alpha)$ has a local minimum. Therefore, the number and location of the $\Pi(\alpha)$ critical points are informative about the number and location of the modes of the mixture density.

Referring back to Section 4, the number of up-down oscillations of $\Pi$ is determined by the *zeroes* of

$$(7) \qquad \Pi'(\alpha) = -\frac{\phi_2''(\alpha)\phi_1'(\alpha) - \phi_1''(\alpha)\phi_2'(\alpha)}{(\phi_2'(\alpha) - \phi_1'(\alpha))^2}.$$

In general, to determine the sign changes of $\Pi'$, we can use any function of $\alpha$ with the same numerator $\phi_2''(\alpha)\phi_1'(\alpha) - \phi_1''(\alpha)\phi_2'(\alpha)$, provided the denominator is a positive function of $\alpha$. In this paper we will use the curvature function $\kappa(\alpha)$ defined by

$$(8) \qquad \kappa(\alpha) = \frac{\phi_2''(\alpha)}{\phi_2(\alpha)}\frac{\phi_1'(\alpha)}{\phi_1(\alpha)} - \frac{\phi_1''(\alpha)}{\phi_1(\alpha)}\frac{\phi_2'(\alpha)}{\phi_2(\alpha)}.$$

We use $\kappa(\alpha)$ as it results in a simple expression for any distribution belonging to the exponential family. We will call it the *curvature function* because its *zeroes* occur at the *zeroes* in curvature of the curve given by $(\phi_1(\alpha), \phi_2(\alpha))$. It is closely related to the mixture curvature measures given by Lindsay [12].

5.2. *Properties of the curvature function $\kappa(\alpha)$.* We now study the role of the curvature function $\kappa(\alpha)$ more closely.

Based on our description of $\Pi$, it is clear that, at the first *zero*, $\alpha_1$, of $\kappa$, the function $\Pi$ has a maximum, at the next $\alpha_2$ a minimum, and so forth. Thus, if we calculate the values of $\Pi(\alpha_j) = \pi_j$, then we can determine from them the ranges of $\pi$ in which we will have modes and antimodes. If the function $\kappa(\alpha)$ does not change its sign in the range $\alpha \in [0, 1]$, then the density is unimodal for all $\pi$. If $\kappa(\alpha)$ has exactly two *zeroes* at $\alpha_1$ and $\alpha_2$, then, for $\pi$ between $\Pi(\alpha_2)$ and $\Pi(\alpha_1)$, the mixture will have two modes, and otherwise there will be just one; and so forth.

In our next result, we develop a simple expression for the curvature. Recall that we defined $S_\alpha = \bar{\alpha}\Sigma_1^{-1} + \alpha\Sigma_2^{-1}$.

THEOREM 3. *Let $g(\mathbf{x})$ be the mixture of two multivariate normal densities. Then*

$$\kappa(\alpha) = [p(\alpha)]^2[1 - \alpha\bar{\alpha}p(\alpha)], \tag{9}$$

*where* $p(\alpha) = (\boldsymbol{\mu}_2 - \boldsymbol{\mu}_1)'\Sigma_1^{-1}S_\alpha^{-1}\Sigma_2^{-1}S_\alpha^{-1}\Sigma_2^{-1}S_\alpha^{-1}\Sigma_1^{-1}(\boldsymbol{\mu}_2 - \boldsymbol{\mu}_1).$

PROOF. Given in the Appendix. □

Now, using the above theorem, since $p(\alpha)$ is always positive, it is clear that the *zeroes* of $\kappa(\alpha)$ are the same as the *zeroes* of $(1 - \alpha\bar{\alpha}p(\alpha))$. For notational ease, let us denote

$$q(\alpha) = 1 - \alpha\bar{\alpha}p(\alpha).$$

By calculation, $q(0) = q(1) = 1$ and, hence, $\kappa$ takes positive values at the two extremes $\alpha = 0$ and $1$. Thus, there are an even number of sign changes of the function $\kappa(\alpha)$ in the range $[0, 1]$, as also indicated by the nature of $\Pi$.

In some special cases one can analytically describe the *zeroes* of $\kappa(\alpha)$. In the corollaries which follow we show that our curvature result can be used to duplicate some of the univariate results found in the literature, and we extend them to certain multivariate situations.

COROLLARY 3. *In the equal variance case $(\Sigma_1 = \Sigma_2)$, $q(\alpha)$ reduces to a quadratic expression*

$$q(\alpha) = 1 - \alpha(1 - \alpha)(\boldsymbol{\mu}_2 - \boldsymbol{\mu}_1)'\Sigma^{-1}(\boldsymbol{\mu}_2 - \boldsymbol{\mu}_1).$$

*$q(\alpha)$ and, hence, $\kappa(\alpha)$ have two roots iff $(\boldsymbol{\mu}_2 - \boldsymbol{\mu}_1)'\Sigma^{-1}(\boldsymbol{\mu}_2 - \boldsymbol{\mu}_1) > 4$, in which case the mixture will be bimodal if and only if $\pi \in (\pi_1, \pi_2)$, where*

$$\frac{1}{\pi_i} = 1 + \frac{\alpha_i}{\bar{\alpha}_i}\frac{\phi_1(\alpha_i)}{\phi_2(\alpha_i)}$$

*and the $\alpha_i$ are the two solutions in $[0, 1]$ of $q(\alpha) = 0$.*

PROOF. Straightforward calculation. □

The above corollary is the generalization of the conditions for bimodality given for the univariate equal variance normal mixture density in [1, 5, 6]. We can also replicate, in the multivariate case, the unequal variance results for univariate densities, provided we make a proportional variance assumption. For this we need the following lemma.

LEMMA 1. *For the proportional variance case* $(\Sigma_2 = \sigma^2 \Sigma_1)$, *the zeroes of* $\kappa(\alpha)$ *in the range* $[0, 1]$ *are the same as the zeroes in* $[0, 1]$ *of the following cubic in* $\alpha$:

$$(10) \qquad q_1(\alpha) = (\sigma^2(1-\alpha) + \alpha)^3 - \alpha(1-\alpha)\mu^2\sigma^2,$$

*where* $\mu^2 = (\boldsymbol{\mu}_2 - \boldsymbol{\mu}_1)'\Sigma_1^{-1}(\boldsymbol{\mu}_2 - \boldsymbol{\mu}_1)$.

PROOF. See the Appendix. □

With the machinery we have now developed, we can offer a simplified proof of the results of Robertson and Fryer [20], while extending those results to higher dimensions.

COROLLARY 4. *Let* $g$ *be the mixture of two normal densities with means* $\boldsymbol{\mu}_1$ *and* $\boldsymbol{\mu}_2$ *and variances* $\Sigma_1$ *and* $\Sigma_2 = \sigma^2\Sigma_1$.

(a) *The density of* $f$ *is unimodal for any mixing proportion* $\pi$ *if*

$$(\boldsymbol{\mu}_2 - \boldsymbol{\mu}_1)'\Sigma_1^{-1}(\boldsymbol{\mu}_2 - \boldsymbol{\mu}_1) \leq \frac{2(1-\sigma^2+\sigma^4)^{3/2} - (2\sigma^6 - 3\sigma^4 - 3\sigma^2 + 2)}{\sigma^2}.$$

(b) *If the parameters do not satisfy the above condition,* $f$ *is bimodal if and only if* $\pi \in (\pi_1, \pi_2)$, *where*

$$\frac{1}{\pi_i} = 1 + \frac{\alpha_i}{\bar{\alpha}_i}\frac{\phi_1(\alpha_i)}{\phi_2(\alpha_i)}$$

*and the* $\alpha_i$ *are the two solutions in* $[0, 1]$ *of*

$$(\sigma^2(1-\alpha) + \alpha)^3 - \alpha(1-\alpha)\sigma^2(\boldsymbol{\mu}_2 - \boldsymbol{\mu}_1)'\Sigma_1^{-1}(\boldsymbol{\mu}_2 - \boldsymbol{\mu}_1) = 0.$$

PROOF. See the Appendix. □

**6. Conclusion.** In this paper we have developed some powerful tools for understanding the topography of a multivariate normal mixture model. The tools are especially powerful in the case $K = 2$, where we can reduce our problem from $D$ dimensions down to one. In any problem one can produce simple plots to investigate the key features of the density. In certain cases we can even describe analytically the number of modes and their locations.

In the process of doing this analysis, we have not discussed how these new results might be used for statistical purposes. We think the possibilities are rich. For example, consider the clustering problem. If we fit a mixture of normals to high-dimensional data, we can associate the components with clusters of data [14]. However, we might also be interested in the information about how well separated two clusters are. A ridgeline elevation plot of their estimated densities will show if they are close enough to each other to form a single mode, in which case we are unlikely to think of them as well separated.

In a model with many components, one might define two components to be linked if they together form a single mode, and then use a map of the linkages to identify how the clusters are associated with each other. (We could also link them if the "mountain passes" are high relative to the peaks.) This could also lead one to a more compact description of the data structure through the construction of supercomponents consisting of linked components, then describing the model as a mixture of a smaller number of supercomponents.

We also note that there are still a number of open mathematical questions. For example, can we find the *zeroes* of the curvature function $\kappa$ analytically in any important special case other than the ones that are given here? Of this we are doubtful. However, finding exact formulae does not seem so important, as we believe it is possible to produce an elementary numerical algorithm to find these points.

In this paper we have reduced the dimension of the mixture modality problem to $\min(K-1, D)$. It was pointed out by a referee that this is exactly the same dimension reduction as occurs in discriminant analysis [2, 7, 8, 9, 19], where one wishes to discriminate between $K$ populations on the basis of $D$ measurements. We think this is an important insight, and that, by using discriminant functions, it will be possible to link our mixture modality study with the area of discriminant analysis in a way that is mutually beneficial.

Through a number of examples in this paper, we have shown that the topography of a multivariate mixture is not like the univariate, as the number of modes can be significantly more than the number of components. As a second question, one might ask if there exists an upper bound for the number of modes, one that can be described as a function of $K$, the number of components, and $D$, the dimension of the multivariate mixture.

Finally, our results are most effective for $K = 2$. It would therefore be useful to establish relationships between the modality structure of the pairs of densities in a mixture and the overall modality of the entire mixture.

NOTE. Datasets and their parameter estimates used in this paper are available for download at **www.stat.psu.edu/~surajit/topography/**.

## APPENDIX: PROOFS

**A.1. Proof of Theorem 1.** Suppose that $\nabla g(\mathbf{x}^*) = 0$, so $\mathbf{x}^*$ is a critical point. Then we have

$$
\begin{aligned}
0 = {}& \pi_1 \phi_1(\mathbf{x}^*) \frac{\nabla \phi_1(\mathbf{x}^*)}{\phi_1(\mathbf{x}^*)} + \pi_2 \phi_2(\mathbf{x}^*) \frac{\nabla \phi_2(\mathbf{x}^*)}{\phi_2(\mathbf{x}^*)} + {} \\
& \cdots + \pi_K \phi_K(\mathbf{x}^*) \frac{\nabla \phi_K(\mathbf{x}^*)}{\phi_K(\mathbf{x}^*)}.
\end{aligned}
\tag{A.1}
$$

If we let

$$(A.2) \qquad \alpha_i = \frac{\pi_i \phi_i(\mathbf{x}^*)}{\pi_1 \phi_1(\mathbf{x}^*) + \pi_2 \phi_2(\mathbf{x}^*) + \cdots + \pi_K \phi_K(\mathbf{x}^*)},$$

then, obviously, $0 \leq \alpha_i \leq 1$ and $\sum_{i=1}^{K} \alpha_i = 1$. Further, note that

$$\frac{\nabla \phi(\mathbf{x}^*; \boldsymbol{\mu}, \Sigma)}{\phi(\mathbf{x}^*; \boldsymbol{\mu}, \Sigma)} = -\Sigma^{-1}(\mathbf{x}^* - \boldsymbol{\mu}).$$

Thus, we have, from equation (A.1), that for every critical value $\mathbf{x}^*$ there exists an $\boldsymbol{\alpha}$ such that

$$(A.3) \quad \alpha_1 \Sigma_1^{-1}(\mathbf{x}^* - \boldsymbol{\mu}_1) + \alpha_2 \Sigma_2^{-1}(\mathbf{x}^* - \boldsymbol{\mu}_2) + \cdots + \alpha_K \Sigma_K^{-1}(\mathbf{x}^* - \boldsymbol{\mu}_K) = 0.$$

Solving this equation for $\mathbf{x}^*$ gives the theorem.

**A.2. Proof of Theorem 2.** First, we note the geometric importance of the vectors $\mathbf{v}_j$. We know that the point $\mathbf{x}(\boldsymbol{\alpha})$ lies in one of the elliptical contours of the density $\phi_j$. At this point the gradient in $\mathbf{x}$ of the density $\phi_j(\mathbf{x})$ is proportional to $\mathbf{v}_j$, and so $\mathbf{v}_j$ is orthogonal to the contour. Thus, if we were to start at $\mathbf{x}(\boldsymbol{\alpha})$ and travel in any direction $\mathbf{w}$ orthogonal to $\mathbf{v}_j$, our path is in the support hyperplane to the elliptically shaped upper set $\{\mathbf{x} : \phi_j(\mathbf{x}) \geq \phi_j(\mathbf{x}^*(\boldsymbol{\alpha}))\}$. Using the fact that the ellipse is convex, our path lies outside the ellipse, and so in the set $\{\mathbf{x} : \phi_j(\mathbf{x}) < \phi_j(\mathbf{x}(\boldsymbol{\alpha}))\}$, except for equality at $\mathbf{x} = \mathbf{x}^*(\boldsymbol{\alpha})$. That is, the point $\mathbf{x}^*(\boldsymbol{\alpha})$ is a local maximum to $\phi_j(\mathbf{x})$ along any path orthogonal to $\mathbf{v}_j$.

Now, suppose that $\mathbf{w}$ satisfies the assumptions of the theorem. It follows from the form of $\mathbf{d}_j$ that

$$(A.4) \qquad \mathbf{w}'(\mathbf{v}_j - \mathbf{v}_K) = 0.$$

However, due to the definition of $\mathbf{x}(\boldsymbol{\alpha})$ we also have $\sum \alpha_j \mathbf{v}_j = 0$, and so $\sum \alpha_j(\mathbf{v}_j - \mathbf{v}_K) = -\mathbf{v}_K$. Putting this together with (A.4) shows that $\mathbf{w}' \mathbf{v}_K = 0$, and so $\mathbf{w}' \mathbf{v}_j = 0$ for $j = 1, \ldots, K$. That is, $\mathbf{w}$ is orthogonal to every $\mathbf{v}_j$, and so, by the first paragraph *every* component of the mixture density $g(x)$ is locally maximized along the given line $\{\mathbf{x}(\boldsymbol{\alpha}) + \delta \mathbf{w} : \delta \in \Re\}$ at $\delta = 0$, and, hence, so is $g(x)$.

**A.3. Proof of Theorem 3.** For simplicity of calculation, we will reuse the notation $\mathbf{v}_j$ introduced in Section 2.2.1. Also, let us denote the first and the second derivatives of a vector $\mathbf{x}$ with respect to the scalar $\alpha$ by $\dot{\mathbf{x}}$ and $\ddot{\mathbf{x}}$, respectively. We will also use the same notation for the first and the second derivatives of the likelihood function.

The equation defining $\mathbf{x}^*$, which is

$$\bar{\alpha} \Sigma_1^{-1}(\mathbf{x}^*(\alpha) - \boldsymbol{\mu}_1) + \alpha \Sigma_2^{-1}(\mathbf{x}^*(\alpha) - \boldsymbol{\mu}_2) = 0,$$

can be written as

$$\text{(A.5)} \qquad \bar{\alpha}\mathbf{v}_1 + \alpha\mathbf{v}_2 = 0$$

$$\text{(A.6)} \qquad \Longrightarrow \quad \mathbf{v}_1 = -\frac{\alpha}{\bar{\alpha}}\mathbf{v}_2.$$

Differentiating (A.5) w.r.t. $\alpha$, we get

$$\text{(A.7)} \qquad \bar{\alpha}\Sigma_1^{-1}\dot{\mathbf{x}}^*(\alpha) - \mathbf{v}_1 + \alpha\Sigma_2^{-1}\dot{\mathbf{x}}^*(\alpha) + \mathbf{v}_2 = 0$$

$$\text{(A.8)} \qquad \Longrightarrow \quad (\bar{\alpha}\Sigma_1^{-1} + \alpha\Sigma_2^{-1})\dot{\mathbf{x}}^*(\alpha) = \mathbf{v}_1 - \mathbf{v}_2$$

$$= \frac{\mathbf{v}_1}{\alpha} \qquad \text{[using (A.6)]}.$$

Taking $l_j(x) = \log(\phi_j(x))$ for $j = 1, 2$, the curvature $\kappa(\alpha)$, given by (8), can be calculated via

$$\frac{\dot{\phi}_j}{\phi_j} = \dot{l}_j \quad \text{and} \quad \frac{\ddot{\phi}_j}{\phi_j} = \ddot{l}_j + [\dot{l}_j]^2.$$

It can be shown that

$$\dot{l}_1 = -\dot{\mathbf{x}}^*(\alpha)\mathbf{v}_1, \qquad\qquad \dot{l}_2 = \frac{\bar{\alpha}}{\alpha}\dot{\mathbf{x}}^*(\alpha)\mathbf{v}_1,$$

$$\ddot{l}_1 = -\ddot{\mathbf{x}}^*(\alpha)\mathbf{v}_1 - \dot{\mathbf{x}}^*(\alpha)\dot{\mathbf{v}}_1, \qquad \ddot{l}_2 = \frac{\bar{\alpha}}{\alpha}\ddot{\mathbf{x}}^*(\alpha)\mathbf{v}_1 + \frac{\bar{\alpha}}{\alpha}\dot{\mathbf{x}}^*(\alpha)\dot{\mathbf{v}}_1 - \dot{\mathbf{x}}^*(\alpha)\mathbf{v}_1\frac{1}{\alpha^2}.$$

After simplification,

$$\kappa(\alpha) = \left(\frac{\dot{\mathbf{x}}^*(\alpha)\mathbf{v}_1}{\alpha}\right)^2\left[1 - \frac{\dot{\mathbf{x}}^*(\alpha)\mathbf{v}_1}{\alpha}\alpha\bar{\alpha}\right].$$

Writing in terms of the original parameters,

$$\frac{\dot{\mathbf{x}}^*(\alpha)\mathbf{v}_1}{\alpha} = (\boldsymbol{\mu}_2 - \boldsymbol{\mu}_1)'\Sigma_1^{-1}S_\alpha^{-1}\Sigma_2^{-1}S_\alpha^{-1}\Sigma_2^{-1}S_\alpha^{-1}\Sigma_1^{-1}(\boldsymbol{\mu}_2 - \boldsymbol{\mu}_1) = p(\alpha).$$

REMARK 6. Note that $p(\alpha)$ function is symmetric in the component labels, as it should be, because one can show by direct calculation that

$$\Sigma_1^{-1}S_\alpha^{-1}\Sigma_2^{-1} = \Sigma_2^{-1}S_\alpha^{-1}\Sigma_1^{-1}.$$

**A.4. Proof of Lemma 1.** $\Sigma_2 = \sigma^2\Sigma_1 \Rightarrow S_\alpha = \Sigma_1^{-1}((1-\alpha) + \frac{\alpha}{\sigma^2})$. Using this relation,

$$q(\alpha) = 1 - \alpha(1-\alpha)\frac{(\boldsymbol{\mu}_2 - \boldsymbol{\mu}_1)'\Sigma_1^{-1}\Sigma_1\Sigma_2^{-1}\Sigma_1\Sigma_2^{-1}\Sigma_1\Sigma_1^{-1}(\boldsymbol{\mu}_2 - \boldsymbol{\mu}_1)}{((1-\alpha) + \alpha/\sigma^2)^3}$$

$$= 1 - \alpha(1-\alpha)\frac{(\boldsymbol{\mu}_2 - \boldsymbol{\mu}_1)'\Sigma_1^{-1}(\boldsymbol{\mu}_2 - \boldsymbol{\mu}_1)}{\sigma^4((1-\alpha) + \alpha/\sigma^2)^3}$$

$$= 1 - \alpha(1-\alpha)\frac{\mu^2\sigma^2}{(\sigma^2(1-\alpha)+\alpha)^3}$$

$$= \frac{1}{(\sigma^2(1-\alpha)+\alpha)^3}q_1(\alpha).$$

Since $(\sigma^2(1-\alpha)+\alpha)^3$ is positive in $[0,1]$, the zeroes of $q(\alpha)$ and $q_1(\alpha)$ are the same in that range.

**A.5. Proof of Corollary 4.** We first prove part (a). Using Lemma 1, we know that the zeroes of $\kappa(\alpha)$ are the same as the zeroes of the cubic $q_1(\alpha) = (\sigma^2(1-\alpha)+\alpha)^3 - \alpha(1-\alpha)\mu^2\sigma^2$. Now, $q_1(\alpha)$ has more than two real zeroes iff the discriminant of $q_1(\alpha)$ is nonnegative. This condition is equivalent to

$$s(\mu) = \mu^4\sigma^2 - \mu^2(-4\sigma^6 + 6\sigma^4 + 6\sigma^2 - 4) - 27\sigma^2(\sigma^2-1)^2 \geq 0,$$

that is,

$$(A.9) \quad s(\mu) = \mu^4\sigma^2 + 2\mu^2(\sigma^2-2)(\sigma^2+1)(2\sigma^2-1) - 27\sigma^2(\sigma^2-1)^2 \geq 0.$$

Now $s(\mu)$ is a quadratic in $\mu^2$ and, among its two zeroes only one is positive. The positive zero is given by

$$\mu_0^2 = \frac{2(1-\sigma^2+\sigma^4)^{3/2} - (2\sigma^6 - 3\sigma^4 - 3\sigma^2 + 2)}{\sigma^2}.$$

Also note that $s(0) < 0$. Thus, $s(\mu)$ is positive when $\mu^2 > \mu_0^2$. Now applying Theorem 3 and Lemma 1, in the proportional variance case we find that $g$ is unimodal if

$$(A.10) \quad \begin{aligned} &(\boldsymbol{\mu}_2 - \boldsymbol{\mu}_1)'\Sigma_1^{-1}(\boldsymbol{\mu}_2 - \boldsymbol{\mu}_1) \\ &\leq \frac{2(1-\sigma^2+\sigma^4)^{3/2} - (2\sigma^6 - 3\sigma^4 - 3\sigma^2 + 2)}{\sigma^2}. \end{aligned}$$

Note that this condition is exactly the condition Robertson and Fryer [20] derived for the univariate case.

Now we prove part (b). We have already noted that the change of curvature of the $\Pi(\alpha)$ occurs at the same place as the solutions to $\kappa(\alpha) = 0$, which in turn are the roots of the cubic equation $q_1(\alpha) = 0$ for the proportional variance case. Let $\alpha_1$ and $\alpha_2$ be the two roots of

$$q_1(\alpha) = (\sigma^2(1-\alpha)+\alpha)^3 - \alpha(1-\alpha)\mu^2\sigma^2 = 0$$

lying between 0 and 1. Thus, the range of $\pi$ for which the density is bimodal can be obtained from the range of $\alpha$ for which $\kappa(\alpha) < 0$, which is the interval $(\alpha_1, \alpha_2)$. Using the relation in (6), the range of $\pi$ can be derived as being the open interval $(\pi_1, \pi_2)$, such that

$$\frac{1}{\pi_i} = 1 + \frac{\alpha_i}{\bar{\alpha}_i}\frac{\phi_1(\alpha_i)}{\phi_2(\alpha_i)} \qquad \text{for } i = 1, 2.$$

Department of Biostatistics
University of North Carolina
Chapel Hill, North Carolina 27599
USA
E-mail: sray@bios.unc.edu

Department of Statistics
Pennsylvania State University
University Park, Pennsylvania 16802
USA
E-mail: bgl@stat.psu.edu